\documentclass[11pt]{amsart}
\usepackage{amsmath,amssymb}
\newtheorem{theorem}{Theorem}

\newtheorem{corollary}[theorem]{Corollary}
\newtheorem{lemma}[theorem]{Lemma}

\begin{document}

\title[The isoperimetric inequality for a minimal submanifold]{The isoperimetric inequality for a minimal submanifold in Euclidean space} 
\author{Simon Brendle}
\address{Department of Mathematics \\ Columbia University \\ New York NY 10027}
\begin{abstract}
We prove a Sobolev inequality which holds on submanifolds in Euclidean space of arbitrary dimension and codimension. This inequality is sharp if the codimension is at most $2$. As a special case, we obtain a sharp isoperimetric inequality for minimal submanifolds in Euclidean space of codimension at most $2$.
\end{abstract}
\thanks{This project was supported by the National Science Foundation under grant DMS-1806190 and by the Simons Foundation.}

\maketitle 

\section{Introduction}

The isoperimetric inequality for a domain in $\mathbb{R}^n$ is one of the most beautiful results in geometry. It has long been conjectured that the isoperimetric inequality still holds if we replace the domain in $\mathbb{R}^n$ by a minimal hypersurface in $\mathbb{R}^{n+1}$. In this paper, we prove this conjecture, as well as a more general inequality which holds for submanifolds of arbitrary dimension and codimension.

\begin{theorem} 
\label{arbitrary.codim}
Let $\Sigma$ be a compact $n$-dimensional submanifold of $\mathbb{R}^{n+m}$ (possibly with boundary $\partial \Sigma$), where $m \geq 2$. Let $f$ be a positive smooth function on $\Sigma$. Then 
\[\int_\Sigma \sqrt{|\nabla^\Sigma f|^2 + f^2 \, |H|^2} + \int_{\partial \Sigma} f \geq n \, \Big ( \frac{(n+m) \, |B^{n+m}|}{m \, |B^m|} \Big )^{\frac{1}{n}} \, \Big ( \int_\Sigma f^{\frac{n}{n-1}} \Big )^{\frac{n-1}{n}}.\] 
Here, $H$ denotes the mean curvature vector of $\Sigma$, and $B^n$ denotes the open unit ball in $\mathbb{R}^n$.
\end{theorem}

Let us consider the special case $m=2$. The standard recursion formula for the volume of the unit ball in Euclidean space gives $(n+2) \, |B^{n+2}| = 2\pi \, |B^n| = 2 \, |B^2| \, |B^n|$. Thus, Theorem \ref{arbitrary.codim} implies a sharp Sobolev inequality for submanifolds of codimension $2$:

\begin{corollary} 
\label{codim.2}
Let $\Sigma$ be a compact $n$-dimensional submanifold of $\mathbb{R}^{n+2}$ (possibly with boundary $\partial \Sigma$), and let $f$ be a positive smooth function on $\Sigma$. Then 
\[\int_\Sigma \sqrt{|\nabla^\Sigma f|^2 + f^2 \, |H|^2} + \int_{\partial \Sigma} f \geq n \, |B^n|^{\frac{1}{n}} \, \Big ( \int_\Sigma f^{\frac{n}{n-1}} \Big )^{\frac{n-1}{n}},\] 
where $H$ denotes the mean curvature vector of $\Sigma$.
\end{corollary}

Finally, we characterize the case of equality in Corollary \ref{codim.2}:

\begin{theorem} 
\label{equality.case}
Let $\Sigma$ be a compact $n$-dimensional submanifold of $\mathbb{R}^{n+2}$ (possibly with boundary $\partial \Sigma$), and let $f$ be a positive smooth function on $\Sigma$. If
\[\int_\Sigma \sqrt{|\nabla^\Sigma f|^2 + f^2 \, |H|^2} + \int_{\partial \Sigma} f = n \, |B^n|^{\frac{1}{n}} \, \Big ( \int_\Sigma f^{\frac{n}{n-1}} \Big )^{\frac{n-1}{n}},\] 
then $f$ is constant and $\Sigma$ is a flat round ball.
\end{theorem}

In particular, if $\Sigma$ is a compact $n$-dimensional minimal submanifold of $\mathbb{R}^{n+2}$, then $\Sigma$ satisfies the sharp isoperimetric inequality 
\[|\partial \Sigma| \geq n \, |B^n|^{\frac{1}{n}} \, |\Sigma|^{\frac{n-1}{n}},\] 
and equality holds if and only if $\Sigma$ is a flat round ball.

Every $n$-dimensional submanifold of $\mathbb{R}^{n+1}$ can be viewed as a submanifold of $\mathbb{R}^{n+2}$. Hence, Corollary \ref{codim.2} and Theorem \ref{equality.case} imply a sharp isoperimetric inequality in codimension $1$.

The isoperimetric inequality on a minimal surface has a long history. In 1921, Torsten Carleman \cite{Carleman} proved that every two-dimensional minimal surface $\Sigma$ which is diffeomorphic to a disk satisfies the sharp isoperimetric inequality $|\partial \Sigma|^2 \geq 4\pi \, |\Sigma|$. Various authors have weakened the topological assumption in Carleman's theorem. In particular, the sharp isoperimetric inequality has been verified for two-dimensional minimal surfaces with connected boundary (see \cite{Hsiung}, \cite{Reid}); for two-dimensional minimal surfaces diffeomorphic to annuli (cf. \cite{Feinberg}, \cite{Osserman-Schiffer}); and for two-dimensional minimal surfaces with two boundary components (cf. \cite{Choe1}, \cite{Li-Schoen-Yau}). On the other hand, using different techniques, Leon Simon showed that every two-dimensional minimal surface satisfies the non-sharp isoperimetric inequality $|\partial \Sigma|^2 \geq 2\pi \, |\Sigma|$ (see \cite{Topping}, Section 4). Stone \cite{Stone} subsequently improved the constant in this inequality: he proved that $|\partial \Sigma|^2 \geq 2\sqrt{2} \, \pi \, |\Sigma|$ for every two-dimensional minimal surface $\Sigma$. We refer to \cite{Choe2} for a survey of these developments. 

In higher dimensions, the famous Michael-Simon Sobolev inequality (cf. \cite{Allard}, Section 7, and \cite{Michael-Simon}) implies an isoperimetric inequality for minimal submanifolds, albeit with a non-sharp constant. Castillon \cite{Castillon} gave an alternative proof of the Michael-Simon Sobolev inequality using methods from optimal transport. Finally, Almgren \cite{Almgren} proved a sharp version of the filling inequality of Federer and Fleming \cite{Federer-Fleming}. In particular, this gives a sharp isoperimetric inequality for area-minimizing submanifolds in all dimensions.

Our method of proof is inspired in part by the Alexandrov-Bakelman-Pucci maximum principle (cf. \cite{Cabre}, \cite{Trudinger}). An alternative way to prove Theorem \ref{arbitrary.codim} would be to use optimal transport; in that case, we would consider the transport map from a thin annulus in $\mathbb{R}^{n+m}$ to the submanifold $\Sigma$ equipped with the measure $f^{\frac{n}{n-1}} \, d\text{\rm vol}$.

\section{Proof of Theorem \ref{arbitrary.codim}}

\label{proof.of.main.result}

Let $\Sigma$ be a compact $n$-dimensional submanifold of $\mathbb{R}^{n+m}$ (possibly with boundary $\partial \Sigma$), where $m \geq 2$. For each point $x \in \Sigma$, we denote by $T_x \Sigma$ and $T_x^\perp \Sigma$ the tangent and normal space to $\Sigma$ at $x$, respectively. Moreover, we denote by $I\!I$ the second fundamental form of $\Sigma$. Recall that $I\!I$ is a symmetric bilinear form on $T_x \Sigma$ which takes values in $T_x^\perp \Sigma$. If $X$ and $Y$ are tangent vector fields on $\Sigma$ and $V$ is a normal vector field along $\Sigma$, then $\langle I\!I(X,Y),V \rangle = \langle \bar{D}_X Y,V \rangle = -\langle \bar{D}_X V,Y \rangle$, where $\bar{D}$ denotes the standard connection on $\mathbb{R}^{n+m}$. The trace of the second fundamental form gives the mean curvature vector, which we denote by $H$. Finally, we denote by $\eta$ the co-normal to $\partial \Sigma$. 

We now turn to the proof of Theorem \ref{arbitrary.codim}. We first consider the special case that $\Sigma$ is connected. By scaling, we may assume that 
\[\int_\Sigma \sqrt{|\nabla^\Sigma f|^2 + f^2 \, |H|^2} + \int_{\partial \Sigma} f = n \int_\Sigma f^{\frac{n}{n-1}}.\] 
Since $\Sigma$ is connected, we can find a function $u: \Sigma \to \mathbb{R}$ with the property that 
\[\text{\rm div}_\Sigma(f \, \nabla^\Sigma u) = n \, f^{\frac{n}{n-1}} - \sqrt{|\nabla^\Sigma f|^2 + f^2 \, |H|^2}\] 
on $\Sigma$ and $\langle \nabla^\Sigma u,\eta \rangle = 1$ at each point on $\partial \Sigma$. Since the function $\sqrt{|\nabla^\Sigma f|^2 + f^2 \, |H|^2}$ is Lipschitz continuous, it follows from standard elliptic regularity theory that the function $u$ is of class $C^{2,\gamma}$ for each $0 < \gamma < 1$ (see \cite{Gilbarg-Trudinger}, Theorem 6.30).

We define 
\begin{align*} 
\Omega &:= \{x \in \Sigma \setminus \partial \Sigma: |\nabla^\Sigma u(x)| < 1\}, \\ 
U &:= \{(x,y): x \in \Sigma \setminus \partial \Sigma, \, y \in T_x^\perp \Sigma, \, |\nabla^\Sigma u(x)|^2 + |y|^2 < 1\}, \\ 
A &:= \{(x,y) \in U: D_\Sigma^2 u(x) - \langle I\!I(x),y \rangle \geq 0\}. 
\end{align*} 
Moreover, we define a map $\Phi: U \to \mathbb{R}^{n+m}$ by 
\[\Phi(x,y) = \nabla^\Sigma u(x) + y\] 
for all $(x,y) \in U$. Note that $\Phi$ is of class $C^{1,\gamma}$ for each $0 < \gamma < 1$. Since $\nabla^\Sigma u(x) \in T_x \Sigma$ and $y \in T_x^\perp \Sigma$ are orthogonal, we obtain $|\Phi(x,y)|^2 = |\nabla^\Sigma u(x)|^2+|y|^2 < 1$ for all $(x,y) \in U$.

\begin{lemma}
\label{Phi.surjective}
The image $\Phi(A)$ is the open unit ball $B^{n+m}$.
\end{lemma}

\textbf{Proof.} 
Clearly, $\Phi(A) \subset \Phi(U) \subset B^{n+m}$. To prove the reverse inclusion, we consider an arbitrary vector $\xi \in \mathbb{R}^{n+m}$ such that $|\xi|<1$. We define a function $w: \Sigma \to \mathbb{R}$ by $w(x) := u(x) - \langle x,\xi \rangle$. Using the Cauchy-Schwarz inequality, we obtain 
\[\langle \nabla^\Sigma w(x),\eta(x) \rangle = \langle \nabla^\Sigma u(x),\eta(x) \rangle - \langle \eta(x),\xi \rangle = 1 - \langle \eta(x),\xi \rangle > 0\] 
for each point $x \in \partial \Sigma$. Consequently, the function $w$ must attain its minimum in the interior of $\Sigma$. Let $\bar{x} \in \Sigma \setminus \partial \Sigma$ be a point in the interior of $\Sigma$ such that $w(\bar{x}) = \inf_{x \in \Sigma} w(x)$. Clearly, $\nabla^\Sigma w(\bar{x}) = 0$. This implies $\xi = \nabla^\Sigma u(\bar{x}) + \bar{y}$ for some $\bar{y} \in T_{\bar{x}}^\perp \Sigma$. Consequently, $|\nabla^\Sigma u(\bar{x})|^2 + |\bar{y}|^2 = |\xi|^2 < 1$. Moreover, we have $D_\Sigma^2 w(\bar{x}) \geq 0$. From this, we deduce that $D_\Sigma^2 u(\bar{x}) - \langle I\!I(\bar{x}),\xi \rangle \geq 0$. Since $\langle I\!I(\bar{x}),\xi \rangle = \langle I\!I(\bar{x}),\nabla^\Sigma u(\bar{x}) + \bar{y} \rangle = \langle I\!I(\bar{x}),\bar{y} \rangle$, we conclude that $D_\Sigma^2 u(\bar{x}) - \langle I\!I(\bar{x}),\bar{y} \rangle \geq 0$. Therefore, $(\bar{x},\bar{y}) \in A$ and $\Phi(\bar{x},\bar{y}) = \xi$. Thus, $B^{n+m} \subset \Phi(A)$. \\

\begin{lemma} 
\label{Jacobian.determinant}
The Jacobian determinant of $\Phi$ is given by 
\[\det D\Phi(x,y) = \det (D_\Sigma^2 u(x) - \langle I\!I(x),y \rangle)\] 
for all $(x,y) \in U$.
\end{lemma}

\textbf{Proof.} 
Fix a point $(\bar{x},\bar{y}) \in U$. Let $\{e_1,\hdots,e_n\}$ be an orthonormal basis of the tangent space $T_{\bar{x}} \Sigma$, and let $(x_1,\hdots,x_n)$ be a local coordinate system on $\Sigma$ such that $\frac{\partial}{\partial x_i} = e_i$ at the point $\bar{x}$. Moreover, let $\{\nu_1,\hdots,\nu_m\}$ denote a local orthonormal frame for the normal bundle $T^\perp \Sigma$. Every normal vector $y$ can be written in the form $y = \sum_{\alpha=1}^m y_\alpha \nu_\alpha$. With this understood, $(x_1,\hdots,x_n,y_1,\hdots,y_m)$ is a local coordinate system on the total space of the normal bundle $T^\perp \Sigma$. We compute 
\begin{align*} 
\Big \langle \frac{\partial \Phi}{\partial x_i}(\bar{x},\bar{y}),e_j \Big \rangle 
&= \langle \bar{D}_{e_i} \nabla^\Sigma u,e_j \rangle + \sum_{\alpha=1}^m \bar{y}_\alpha \, \langle \bar{D}_{e_i} \nu_\alpha,e_j \rangle \\ 
&= (D_\Sigma^2 u)(e_i,e_j) - \langle I\!I(e_i,e_j),\bar{y} \rangle. 
\end{align*} 
In the last step, we have used the identity $\langle I\!I(e_i,e_j),\nu_\alpha \rangle = -\langle \bar{D}_{e_i} \nu_\alpha,e_j \rangle$. Moreover, 
\[\Big \langle \frac{\partial \Phi}{\partial y_\alpha}(\bar{x},\bar{y}),e_j \Big \rangle = \langle \nu_\alpha,e_j \rangle = 0\] 
and  
\[\Big \langle \frac{\partial \Phi}{\partial y_\alpha}(\bar{x},\bar{y}),\nu_\beta \Big \rangle = \langle \nu_\alpha,\nu_\beta \rangle = \delta_{\alpha\beta}.\] 
Thus, we conclude that 
\[\det D\Phi(\bar{x},\bar{y}) = \det \begin{bmatrix} D_\Sigma^2 u(\bar{x}) - \langle I\!I(\bar{x}),\bar{y} \rangle & 0 \\ * & \text{\rm id} \end{bmatrix} = \det (D_\Sigma^2 u(\bar{x}) - \langle I\!I(\bar{x}),\bar{y} \rangle).\] 
This proves the assertion. \\

\begin{lemma} 
\label{Jacobian.estimate}
The Jacobian determinant of $\Phi$ satisfies 
\[0 \leq \det D\Phi(x,y) \leq f(x)^{\frac{n}{n-1}}\] 
for all $(x,y) \in A$.
\end{lemma}

\textbf{Proof.} 
Consider a point $(x,y) \in A$. Using the inequality $|\nabla^\Sigma u(x)|^2+|y|^2 < 1$ and the Cauchy-Schwarz inequality, we obtain 
\begin{align*} 
&-\langle \nabla^\Sigma f(x),\nabla^\Sigma u(x) \rangle - f(x) \, \langle H(x),y \rangle \\ 
&\leq \sqrt{|\nabla^\Sigma f(x)|^2 + f(x)^2 \, |H(x)|^2} \, \sqrt{|\nabla^\Sigma u(x)|^2+|y|^2} \\ 
&\leq \sqrt{|\nabla^\Sigma f(x)|^2 + f(x)^2 \, |H(x)|^2}. 
\end{align*}
Using the identity $\text{\rm div}_\Sigma(f \, \nabla^\Sigma u) = n \, f^{\frac{n}{n-1}} - \sqrt{|\nabla^\Sigma f|^2 + f^2 \, |H|^2}$, we deduce that 
\begin{align*} 
&\Delta_\Sigma u(x) - \langle H(x),y \rangle \\ 
&= n \, f(x)^{\frac{1}{n-1}} - f(x)^{-1} \, \sqrt{|\nabla^\Sigma f(x)|^2 + f(x)^2 \, |H(x)|^2} \\ 
&- f(x)^{-1} \, \langle \nabla^\Sigma f(x),\nabla^\Sigma u(x) \rangle - \langle H(x),y \rangle \\ 
&\leq n \, f(x)^{\frac{1}{n-1}}. 
\end{align*} 
Moreover, $D_\Sigma^2 u(x) - \langle I\!I(x),y \rangle \geq 0$ since $(x,y) \in A$. Hence, the arithmetic-geometric mean inequality implies
\[0 \leq \det (D_\Sigma^2 u(x) - \langle I\!I(x),y \rangle) \leq \Big ( \frac{\text{\rm tr}(D_\Sigma^2 u(x) - \langle I\!I(x),y \rangle)}{n} \Big )^n \leq f(x)^{\frac{n}{n-1}}.\] 
Using Lemma \ref{Jacobian.determinant}, we conclude that $0 \leq \det D\Phi(x,y) \leq f(x)^{\frac{n}{n-1}}$. This completes the proof of Lemma \ref{Jacobian.estimate}. \\ 

We now continue with the proof of Theorem \ref{arbitrary.codim}. Using Lemma \ref{Phi.surjective} and Lemma \ref{Jacobian.estimate}, we obtain 
\begin{align*} 
&|B^{n+m}| \, (1-\sigma^{n+m}) \\ 
&= \int_{\{\xi \in \mathbb{R}^{n+m}: \sigma^2 < |\xi|^2 < 1\}} 1 \, d\xi \\ 
&\leq \int_\Omega \bigg ( \int_{\{y \in T_x^\perp \Sigma: \sigma^2 < |\Phi(x,y)|^2 < 1\}} |\det D\Phi(x,y)| \, 1_A(x,y) \, dy \bigg ) \, d\text{\rm vol}(x) \\ 
&\leq \int_\Omega \bigg ( \int_{\{y \in T_x^\perp \Sigma: \sigma^2 < |\nabla^\Sigma u(x)|^2+|y|^2 < 1\}} f(x)^{\frac{n}{n-1}} \, dy \bigg ) \, d\text{\rm vol}(x) \\ 
&= |B^m| \int_\Omega \Big [ (1-|\nabla^\Sigma u(x)|^2)^{\frac{m}{2}} -  (\sigma^2-|\nabla^\Sigma u(x)|^2)_+^{\frac{m}{2}} \Big ] \, f(x)^{\frac{n}{n-1}} \, d\text{\rm vol}(x)
\end{align*} 
for all $0 \leq \sigma < 1$. Since $m \geq 2$, the mean value theorem gives $b^{\frac{m}{2}} - a^{\frac{m}{2}} \leq \frac{m}{2} \, (b-a)$ for $0 \leq a \leq b \leq 1$. Consequently,  
\begin{align*} 
&(1-|\nabla^\Sigma u(x)|^2)^{\frac{m}{2}} -  (\sigma^2-|\nabla^\Sigma u(x)|^2)_+^{\frac{m}{2}} \\ 
&\leq \frac{m}{2} \, \Big [ (1-|\nabla^\Sigma u(x)|^2) -  (\sigma^2-|\nabla^\Sigma u(x)|^2)_+ \Big ] \leq \frac{m}{2} \, (1-\sigma^2) 
\end{align*}
for all $x \in \Omega$ and all $0 \leq \sigma < 1$. Putting these facts, together, we obtain
\[|B^{n+m}| \, (1-\sigma^{n+m}) \leq \frac{m}{2} \, |B^m| \, (1-\sigma^2) \int_\Omega f^{\frac{n}{n-1}}\] 
for all $0 \leq \sigma < 1$. In the next step, we divide by $1-\sigma$ and take the limit as $\sigma \to 1$. This gives 
\[(n+m) \, |B^{n+m}| \leq m \, |B^m| \int_\Omega f^{\frac{n}{n-1}} \leq m \, |B^m| \int_\Sigma f^{\frac{n}{n-1}}.\] 
On the other hand, $\int_\Sigma \sqrt{|\nabla^\Sigma f|^2 + f^2 \, |H|^2} + \int_{\partial \Sigma} f = n \int_\Sigma f^{\frac{n}{n-1}}$ in view of our normalization. Thus, we conclude that 
\begin{align*} 
&\int_\Sigma \sqrt{|\nabla^\Sigma f|^2 + f^2 \, |H|^2} + \int_{\partial \Sigma} f \\ 
&= n \int_\Sigma f^{\frac{n}{n-1}} \geq n \, \Big ( \frac{(n+m) \, |B^{n+m}|}{m \, |B^m|} \Big )^{\frac{1}{n}} \, \Big ( \int_\Sigma f^{\frac{n}{n-1}} \Big )^{\frac{n-1}{n}}. 
\end{align*} 
This proves Theorem \ref{arbitrary.codim} in the special case when $\Sigma$ is connected. 

It remains to consider the case when $\Sigma$ is disconnected. In that case, we apply the inequality to each individual connected component of $\Sigma$, and take the sum over all connected components. Since 
\[a^{\frac{n-1}{n}} + b^{\frac{n-1}{n}} > a \, (a+b)^{-\frac{1}{n}} + b \, (a+b)^{-\frac{1}{n}} = (a+b)^{\frac{n-1}{n}}\] 
for $a,b>0$, we conclude that 
\[\int_\Sigma \sqrt{|\nabla^\Sigma f|^2 + f^2 \, |H|^2} + \int_{\partial \Sigma} f > n \, \Big ( \frac{(n+m) \, |B^{n+m}|}{m \, |B^m|} \Big )^{\frac{1}{n}} \, \Big ( \int_\Sigma f^{\frac{n}{n-1}} \Big )^{\frac{n-1}{n}}.\] 
if $\Sigma$ is disconnected. This completes the proof of Theorem \ref{arbitrary.codim}. \\

\section{Proof of Theorem \ref{equality.case}}

Suppose that $\Sigma$ is a compact $n$-dimensional submanifold in $\mathbb{R}^{n+2}$ (possibly with boundary $\partial \Sigma$), and $f$ is a positive smooth function on $\Sigma$ satisfying 
\[\int_\Sigma \sqrt{|\nabla^\Sigma f|^2 + f^2 \, |H|^2} + \int_{\partial \Sigma} f = n \, |B^n|^{\frac{1}{n}} \, \Big ( \int_\Sigma f^{\frac{n}{n-1}} \Big )^{\frac{n-1}{n}}.\] 
Clearly, $\Sigma$ must be connected. 

By scaling, we may arrange that $\int_\Sigma \sqrt{|\nabla^\Sigma f|^2 + f^2 \, |H|^2} + \int_{\partial \Sigma} f = n \, |B^n|$ and $\int_\Sigma f^{\frac{n}{n-1}} = |B^n|$. In particular, 
\[\int_\Sigma \sqrt{|\nabla^\Sigma f|^2 + f^2 \, |H|^2} + \int_{\partial \Sigma} f = n \int_\Sigma f^{\frac{n}{n-1}}.\] 
Let $u: \Sigma \to \mathbb{R}$ denote the solution of the equation 
\[\text{\rm div}_\Sigma(f \, \nabla^\Sigma u) = n \, f^{\frac{n}{n-1}} - \sqrt{|\nabla^\Sigma f|^2 + f^2 \, |H|^2}\] 
on $\Sigma$ with boundary condition $\langle \nabla^\Sigma u,\eta \rangle = 1$ on $\partial \Sigma$. Let $\Omega$, $U$, $A$, and $\Phi: U \to \mathbb{R}^{n+2}$ be defined as in Section \ref{proof.of.main.result}.

\begin{lemma}
\label{W}
Suppose that $\bar{x} \in \Omega$, $\bar{y} \in T_{\bar{x}}^\perp \Sigma$, $|\nabla^\Sigma u(\bar{x})|^2 + |\bar{y}|^2 = 1$, and $D_\Sigma^2 u(\bar{x}) - \langle I\!I(\bar{x}),\bar{y} \rangle \neq f(\bar{x})^{\frac{1}{n-1}} \, g$. Then there exists a real number $\varepsilon \in (0,1)$ and an open neighborhood $W$ of the point $(\bar{x},\bar{y})$ such that $\det D\Phi(x,y) \leq (1-\varepsilon) \, f(x)^{\frac{n}{n-1}}$ for all $(x,y) \in A \cap W$. 
\end{lemma}

\textbf{Proof.} 
We distinguish two cases: 

\textit{Case 1:} Suppose that $D_\Sigma^2 u(\bar{x}) - \langle I\!I(\bar{x}),\bar{y} \rangle \geq 0$. Since $|\nabla^\Sigma u(\bar{x})|^2 + |\bar{y}|^2 = 1$, the Cauchy-Schwarz inequality implies 
\[-\langle \nabla^\Sigma f(\bar{x}),\nabla^\Sigma u(\bar{x}) \rangle - f(\bar{x}) \, \langle H(\bar{x}),\bar{y} \rangle \leq \sqrt{|\nabla^\Sigma f(\bar{x})|^2 + f(\bar{x})^2 \, |H(\bar{x})|^2}.\] 
Using the identity $\text{\rm div}_\Sigma(f \, \nabla^\Sigma u) = n \, f^{\frac{n}{n-1}} - \sqrt{|\nabla^\Sigma f|^2 + f^2 \, |H|^2}$, we obtain 
\[\Delta_\Sigma u(\bar{x}) - \langle H(\bar{x}),\bar{y} \rangle \leq n \, f(\bar{x})^{\frac{1}{n-1}}.\] 
Since $D_\Sigma^2 u(\bar{x}) - \langle I\!I(\bar{x}),\bar{y} \rangle \geq 0$ and $D_\Sigma^2 u(\bar{x}) - \langle I\!I(\bar{x}),\bar{y} \rangle \neq f(\bar{x})^{\frac{1}{n-1}} \, g$, the arithmetic-geometric mean inequality gives 
\[\det (D_\Sigma^2 u(\bar{x}) - \langle I\!I(\bar{x}),\bar{y} \rangle) < f(\bar{x})^{\frac{n}{n-1}}.\] 
Let us choose a real number $\varepsilon \in (0,1)$ such that $\det (D_\Sigma^2 u(\bar{x}) - \langle I\!I(\bar{x}),\bar{y} \rangle) < (1-\varepsilon) \, f(\bar{x})^{\frac{n}{n-1}}$. Since $u$ is of class $C^{2,\gamma}$, we can find an open neighborhood $W$ of $(\bar{x},\bar{y})$ such that $\det (D_\Sigma^2 u(x) - \langle I\!I(x),y \rangle) \leq (1-\varepsilon) \, f(x)^{\frac{n}{n-1}}$ for all $(x,y) \in W$. Using Lemma \ref{Jacobian.determinant}, we obtain $\det D\Phi(x,y) \leq (1-\varepsilon) \, f(x)^{\frac{n}{n-1}}$ for all $(x,y) \in U \cap W$.

\textit{Case 2:} Suppose that the smallest eigenvalue of $D_\Sigma^2 u(\bar{x}) - \langle I\!I(\bar{x}),\bar{y} \rangle$ is strictly negative. Since $u$ is of class $C^{2,\gamma}$, we can find an open neighborhood $W$ of $(\bar{x},\bar{y})$ with the property that the smallest eigenvalue of $D_\Sigma^2 u(x) - \langle I\!I(x),y \rangle$ is strictly negative for all $(x,y) \in W$. Consequently, $A \cap W = \emptyset$. This completes the proof of Lemma \ref{W}. \\

\begin{lemma}
\label{equality.1}
We have $D_\Sigma^2 u(x) - \langle I\!I(x),y \rangle = f(x)^{\frac{1}{n-1}} \, g$ for all $x \in \Omega$ and all $y \in T_x^\perp \Sigma$ satisfying $|\nabla^\Sigma u(x)|^2 + |y|^2 = 1$.
\end{lemma}

\textbf{Proof.} 
We argue by contradiction. Suppose that there exists a point $\bar{x} \in \Omega$ and a vector $\bar{y} \in T_{\bar{x}}^\perp \Sigma$ such that $|\nabla^\Sigma u(\bar{x})|^2 + |\bar{y}|^2 = 1$ and $D_\Sigma^2 u(\bar{x}) - \langle I\!I(\bar{x}),\bar{y} \rangle \neq f(\bar{x})^{\frac{1}{n-1}} \, g$. By Lemma \ref{W}, we can find a real number $\varepsilon \in (0,1)$ and an open neighborhood $W$ of the point $(\bar{x},\bar{y})$ such that $\det D\Phi(x,y) \leq (1-\varepsilon) \, f(x)^{\frac{n}{n-1}}$ for all $(x,y) \in A \cap W$. Using Lemma \ref{Jacobian.estimate}, we deduce that 
\[0 \leq \det D\Phi(x,y) \leq (1 - \varepsilon \cdot 1_W(x,y)) \, f(x)^{\frac{n}{n-1}}\] 
for all $(x,y) \in A$. Arguing as in Section \ref{proof.of.main.result}, we obtain 
\begin{align*} 
&|B^{n+2}| \, (1-\sigma^{n+2}) \\ 
&= \int_{\{\xi \in \mathbb{R}^{n+2}: \sigma^2 < |\xi|^2 < 1\}} 1 \, d\xi \\ 
&\leq \int_\Omega \bigg ( \int_{\{y \in T_x^\perp \Sigma: \sigma^2 < |\Phi(x,y)|^2 < 1\}} |\det D\Phi(x,y)| \, 1_A(x,y) \, dy \bigg ) \, d\text{\rm vol}(x) \\ 
&\leq \int_\Omega \bigg ( \int_{\{y \in T_x^\perp \Sigma: \sigma^2 < |\nabla^\Sigma u(x)|^2+|y|^2 < 1\}} (1 - \varepsilon \cdot 1_W(x,y)) \, f(x)^{\frac{n}{n-1}} \, dy \bigg ) \, d\text{\rm vol}(x) \\ 
&= |B^2| \int_\Omega \Big [ (1-|\nabla^\Sigma u(x)|^2) -  (\sigma^2-|\nabla^\Sigma u(x)|^2)_+ \Big ] \, f(x)^{\frac{n}{n-1}} \, d\text{\rm vol}(x) \\ 
&- \varepsilon \int_\Omega \bigg ( \int_{\{y \in T_x^\perp \Sigma: \sigma^2 < |\nabla^\Sigma u(x)|^2+|y|^2 < 1\}} 1_W(x,y) \, f(x)^{\frac{n}{n-1}} \, dy \bigg ) \, d\text{\rm vol}(x) \\ 
&\leq |B^2| \, (1-\sigma^2) \int_\Omega f(x)^{\frac{n}{n-1}} \, d\text{\rm vol}(x) \\ 
&- \varepsilon \int_\Omega \bigg ( \int_{\{y \in T_x^\perp \Sigma: \sigma^2 < |\nabla^\Sigma u(x)|^2+|y|^2 < 1\}} 1_W(x,y) \, f(x)^{\frac{n}{n-1}} \, dy \bigg ) \, d\text{\rm vol}(x) 
\end{align*} 
for all $0 \leq \sigma < 1$. Dividing by $1-\sigma$ and taking the limit as $\sigma \to 1$ gives 
\[(n+2) \, |B^{n+2}| < 2 \, |B^2| \int_\Omega f^{\frac{n}{n-1}} \leq 2 \, |B^2| \int_\Sigma f^{\frac{n}{n-1}} = 2 \, |B^2| \, |B^n|.\] 
This contradicts the fact that $(n+2) \, |B^{n+2}| = 2 \, |B^2| \, |B^n|$. \\

\begin{lemma}
\label{equality.2}
We have $D_\Sigma^2 u(x) = f(x)^{\frac{1}{n-1}} \, g$ and $I\!I(x) = 0$ for all $x \in \Omega$. 
\end{lemma}

\textbf{Proof.} 
Lemma \ref{equality.1} implies $D_\Sigma^2 u(x) - \langle I\!I(x),y \rangle = f(x)^{\frac{1}{n-1}} \, g$ for all $x \in \Omega$ and all $y \in T_x^\perp \Sigma$ satisfying $|\nabla^\Sigma u(x)|^2 + |y|^2 = 1$. Replacing $y$ by $-y$ gives $D_\Sigma^2 u(x) + \langle I\!I(x),y \rangle = f(x)^{\frac{1}{n-1}} \, g$ for all $x \in \Omega$ and all $y \in T_x^\perp \Sigma$ satisfying $|\nabla^\Sigma u(x)|^2 + |y|^2 = 1$. Consequently, $D_\Sigma^2 u(x) = f(x)^{\frac{1}{n-1}} \, g$ and $\langle I\!I(x),y \rangle = 0$ for all $x \in \Omega$ and all $y \in T_x^\perp \Sigma$ satisfying $|\nabla^\Sigma u(x)|^2 + |y|^2 = 1$. From this, the assertion follows. \\

\begin{lemma}
\label{equality.3}
We have $\nabla^\Sigma f(x) = 0$ for all $x \in \Omega$.
\end{lemma}

\textbf{Proof.} 
Using Lemma \ref{equality.2}, we obtain $\Delta_\Sigma u = n \, f^{\frac{1}{n-1}}$ at each point in $\Omega$. This implies $\text{\rm div}_\Sigma(f \, \nabla^\Sigma u) = n \, f^{\frac{n}{n-1}} + \langle \nabla^\Sigma f,\nabla^\Sigma u \rangle$ at each point in $\Omega$. On the other hand, by definition of $u$, we have $\text{\rm div}_\Sigma(f \, \nabla^\Sigma u) = n \, f^{\frac{n}{n-1}} - |\nabla^\Sigma f|$ at each point in $\Omega$. Consequently, $\langle \nabla^\Sigma f,\nabla^\Sigma u \rangle = -|\nabla^\Sigma f|$ at each point in $\Omega$. Since $|\nabla^\Sigma u| < 1$ at each point in $\Omega$, we conclude that $\nabla^\Sigma f = 0$ at each point in $\Omega$. \\

\begin{lemma}
\label{density}
The set $\Omega$ is dense in $\Sigma$.
\end{lemma}

\textbf{Proof.} 
We argue by contradiction. Suppose that $\Omega$ is not dense in $\Sigma$. Then $\int_\Omega f^{\frac{n}{n-1}} < \int_\Sigma f^{\frac{n}{n-1}}$. Hence, the arguments in Section \ref{proof.of.main.result} imply 
\[(n+2) \, |B^{n+2}| \leq 2 \, |B^2| \int_\Omega f^{\frac{n}{n-1}} < 2 \, |B^2| \int_\Sigma f^{\frac{n}{n-1}} = 2 \, |B^2| \, |B^n|.\] 
This contradicts the fact that $(n+2) \, |B^{n+2}| = 2 \, |B^2| \, |B^n|$. \\

Using Lemma \ref{equality.2}, Lemma \ref{equality.3}, and Lemma \ref{density}, we conclude that $D_\Sigma^2 u = f^{\frac{1}{n-1}} \, g$, $I\!I = 0$, and $\nabla^\Sigma f = 0$ at each point on $\Sigma$. Since $\Sigma$ is connected and $\nabla^\Sigma f = 0$ at each point on $\Sigma$, it follows that $f = \lambda^{n-1}$ for some positive constant $\lambda$. Since $\Sigma$ is connected and $I\!I = 0$ at each point on $\Sigma$, $\Sigma$ is contained in an $n$-dimensional plane $P$. Since $D_\Sigma^2 u = f^{\frac{1}{n-1}} \, g = \lambda \, g$ at each point on $\Sigma$, the function $u$ must be of the form $u(x) = \frac{1}{2} \, \lambda \, |x-p|^2 + c$ for some point $p \in P$ and some constant $c$. On the other hand, we know that $|\nabla^\Sigma u| < 1$ at each point on $\Omega$. Using Lemma \ref{density}, it follows that $|\nabla^\Sigma u| \leq 1$ at each point on $\Sigma$. This implies $\Sigma \subset \{x \in P: \lambda \, |x-p| \leq 1\}$. Since $\lambda^n \, |\Sigma| = \int_\Sigma f^{\frac{n}{n-1}} = |B^n|$, we conclude that $\Sigma = \{x \in P: \lambda \, |x-p| \leq 1\}$. This completes the proof of Theorem \ref{equality.case}.

\end{document}